\documentclass[a4paper,11pt]{amsart}

\usepackage{amsmath, amsthm, amssymb, graphicx, subfigure, amsrefs}

\title{Divided Differences of Implicit Functions}
\author{Georg Muntingh}
\address{CMA / Department of Mathematics, University of Oslo, P.O. Box 1053, Blindern, N-0316, Oslo, Norway}
\email{georgmu@math.uio.no}
\author{Michael Floater}
\address{CMA / Department of Informatics, University of Oslo, P.O. Box 1053, Blindern, N-0316, Oslo, Norway}
\email{michaelf@ifi.uio.no}

\date{\today}

\def\fg#1#2{\frac{[#1#2;#2]g}{[#1;#1#2]g}}
\def\PPP{\mathcal{P}}
\def\ff{\mathbf{f}}
\def\ii{\mathbf{i}}
\def\RR{\mathbb{R}}

\newtheorem{theorem}{Theorem}
\theoremstyle{definition}
\newtheorem{example}{Example}

\begin{document}

\begin{abstract} Under general conditions, the equation $g(x,y) = 0$ implicitly defines $y$
locally as a function of $x$.
In this article, we express divided differences of
$y$ in terms of bivariate divided differences of $g$,
generalizing a recent result on divided differences of inverse functions.
\end{abstract}

\maketitle

\section{Introduction}\label{sec:Introduction}
Divided differences can be viewed as a discrete analogue of
derivatives and are commonly used in approximation theory,
see \cite{Boor} for a survey.

Recently, the second author and Lyche established two univariate chain rules for divided differences \cite{DividedDiffChain}, both of which can be viewed as analogous to Fa\`a di Bruno's formula for differentiating composite functions \cite{FaaDiBruno,Johnson-CuriousHistory}. One of these formulas was simultaneously discovered by Wang and Xu \cite{WangXu}. In a follow-up preprint, the other chain rule was generalized to the composition of vector-valued functions of several variables \cite{ChainRuleMultivariate}, yielding a formula analogous to a multivariate version of Fa\`a di Bruno's formula \cite{MultivariateFaaDiBruno}.

In \cite{DividedDiffInverse}, the univariate chain rule was applied
to find a formula for divided differences of the inverse of a function.
In Theorem \ref{thm:main}, the Main Theorem of this paper, we use the multivariate chain rule
to prove a similar formula for divided differences of
implicitly defined functions. Equation \ref{eq:inverse} shows that the
formula for divided differences of inverse functions in
\cite{DividedDiffInverse} follows as a special case.

More precisely, let $y$ be a function that is defined implicitly by a
function $g:\RR^2\to \RR$ via $g\big(x,y(x)\big)=0$ and $\frac{\partial g}{\partial y}\big(x,y(x)\big)\neq 0$,
for every $x$ in an open interval $U\subset \RR$.
Then the Main Theorem states that for any
\[ x_0,\ldots,x_n \in U,
\qquad y_0 := y(x_0), \ldots, y_n := y(x_n) \in y(U) \]
we can express $[x_0,\ldots,x_n]y$ as a sum of terms involving the divided differences
$[x_{i_0},\ldots,x_{i_s};y_{i_s},\ldots y_{i_r}]g$,
with $0\leq i_0 < i_1 < \cdots < i_r \leq n$.

In Section \ref{sec:DividedDifferences}, we define these divided differences and explain our notation. In Section \ref{sec:RecursiveFormula}, we apply the multivariate chain rule to derive a formula that recursively expresses divided differences of $y$ in terms of divided differences of $g$ and lower order divided differences of $y$. Finally, in Section \ref{sec:formula}, we solve this recursive formula to obtain a formula that expresses divided differences of $y$ solely in terms of divided differences of $g$. We end the section with applying the Main Theorem in some special cases. 

\section{Divided Differences}\label{sec:DividedDifferences}
Let $[x_0,\ldots,x_n]f$ denote the \emph{divided difference} of a function
$f:(a, b)\to \RR$ at the points $x_0,\ldots,x_n$,
with $a < x_0 \leq \cdots \leq x_n < b$.
If all inequalities are strict, this notion is recursively defined by
$[x_0]f := f(x_0)$ and
\[ [x_0,\ldots, x_n]f =
  \frac{[x_1,\ldots, x_n]f - [x_0,\ldots, x_{n-1}]f}{x_n - x_0}
  \qquad \textup{if}~n>0. \]
If some of the $\{x_i\}$ coincide, we define $[x_0,\ldots, x_n]f$
as the limit of this formula when the distances between these $\{x_i\}$
become arbitrary small, provided $f$ is sufficiently smooth there.
In particular, when $x_0 = \cdots = x_n$,
one can show that $[x_0,\ldots,x_n]f = f^{(n)}(x_0) / n!$\,.
For given $i_0,\ldots,i_k$ satisfying
$i_0\leq i_1\leq \cdots \leq i_k$, we shall sometimes shorten notation to
\begin{equation}\label{eq:divdiff-short1}
[i_0i_1\cdots i_k]f := [x_{i_0},x_{i_1},\ldots x_{i_k}]f.
\end{equation}

The above definitions generalize to bivariate divided differences as follows.
Let $f: U\to \RR$ be defined on some $2$-dimensional interval
\[ U = (a_1,b_1)\times (a_2,b_2)\subset \RR^2.\]
Suppose we are given $m, n \ge 0$ and
points $x_0, \ldots, x_m \in (a_1, b_1)$
satisfying $x_0 < \cdots < x_m$
and $y_0, \ldots, y_m \in (a_2, b_2)$
satisfying $y_0 < \cdots < y_m$.
The Cartesian product 
\[ \{ x_0, \ldots, x_m \} \times \{ y_0, \ldots, y_n \} \]
defines a rectangular grid of points in $U$.
The \emph{(bivariate) divided difference} of $f$ at this grid,
denoted by
\begin{equation}\label{eq:bivdd}
 [ x_0, \ldots, x_m; y_0, \ldots, y_n ]f,
\end{equation}
can be defined recursively as follows.
If $m = n = 0$, the grid consists of only one point
$(x_0,y_0)$, and we define
$[x_0;y_0]f := f(x_0,y_0)$
as the value of $f$ at this point.
In case $m > 0$, we can define (\ref{eq:bivdd}) as
\[ \frac{[ x_1, \ldots, x_m; y_0, \ldots, y_n ]f - [ x_0, \ldots, x_{m-1}; y_0, \ldots, y_n ]f}{x_m - x_0}, \]
or if $n > 0$, as
\[ \frac{[ x_0, \ldots, x_m; y_1, \ldots, y_n ]f - [ x_0, \ldots, x_m; y_0, \ldots, y_{n-1} ]f}{y_n - y_0}. \]
If both $m>0$ and $n>0$
the divided difference (\ref{eq:bivdd}) is uniquely defined by either
recursion formula.

As for univariate divided differences, we can let some of the points coalesce
by taking limits, as long as $f$ is sufficiently smooth.
In particular when $x_0 = \cdots = x_m$ and $y_0 = \cdots = y_n$,
this legitimates the notation
\[ [x_0, \ldots, x_m; y_0, \ldots, y_n]f 
  := \frac{1}{m!n!} \frac{\partial^{m+n} f}
      {\partial x^m \partial y^n} (x_0,y_0). \]

Similarly to Equation \ref{eq:divdiff-short1}, we shall more often than not shorten the notation for bivariate divided differences to
\begin{equation}\label{eq:divdiff-short2}
[ i_0i_1\cdots i_s; j_0 j_1 \cdots j_t]f := [x_{i_0},x_{i_1},\ldots,x_{i_s};y_{j_0},y_{j_1},\ldots,y_{j_t}]f.
\end{equation}

\section{A Recursive Formula for Implicit Functions}\label{sec:RecursiveFormula}
Let $y$ be a function implicitly defined by $g\big(x,y(x)\big) = 0$ as in Section \ref{sec:Introduction}.
The first step in expressing divided differences of $y$ in terms of those
of $g$ is to express those of $g$ in terms of those of $y$.
This link is provided by a special case of the
the multivariate chain rule of \cite{ChainRuleMultivariate}.
Let
$\RR\stackrel{\ff}{\longrightarrow}\RR^2\stackrel{g}{\longrightarrow}\RR$
be a composition of sufficiently smooth functions
$\ff=(\phi,\psi)$ and $g$.
In this case, the formula of \cite{ChainRuleMultivariate} for $n \geq 1$ is
\begin{equation}\label{eq:SMultivariateChainRule}
  [x_0,x_1,\ldots,x_n](g\circ \ff) =  
   \sum_{k=1}^n\ \sum_{0=i_0<i_1<\cdots <i_k = n}\ \sum_{s=0}^k
\end{equation}
\[ \qquad [\phi(x_{i_0}),\phi(x_{i_1}), \ldots,\phi(x_{i_s});
           \psi(x_{i_s}),\psi(x_{i_{s+1}}), \ldots,\psi(x_{i_k})]g \]
\[ \qquad \times \prod_{l=1  }^s [x_{i_{l-1}},x_{i_{l-1} + 1}, \ldots,x_{i_l}]\phi
          \prod_{l=s+1}^k [x_{i_{l-1}},x_{i_{l-1} + 1}, \ldots,x_{i_l}]\psi. \]
Now we choose $\ff$ to be the graph of a function $y$,
i.e., $\ff: x\mapsto \big(\phi(x),\psi(x)\big) = \big(x, y(x)\big)$.
Then the divided differences of $\phi$ of order greater than one are zero,
implying that the summand is zero unless
$(i_0,i_1,\ldots,i_s) = (0,1,\ldots,s)$; below this condition is realized by restricting the third sum in Equation \ref{eq:SMultivariateChainRule} to integers $s$ that satisfy $s = i_s - i_0$.
Since additionally divided differences of $\phi$ of order one are one, we obtain
\begin{equation}\label{eq:CR1}
[x_0,x_1,\ldots,x_n]g\big(\cdot,y(\cdot)\big) =  \sum_{k=1}^n\ \sum_{0=i_0<i_1<\cdots <i_k = n}\ \sum_{\substack{s = 0 \\s = i_s - i_0}}^k
\end{equation}
\[ \hfill [x_0,x_1,\ldots,x_s; y_{i_s},y_{i_{s+1}},\ldots,y_{i_k}]g 
   \prod_{l=s+1}^k [x_{i_{l-1}},x_{i_{l-1}+1},\ldots,x_{i_l}]y, \]
where $y_j := y(x_j)$ for $j = 0,1,\ldots,n$. For example, when $n=1$ this formula becomes
\[ [x_0,x_1]g\big(\cdot,y(\cdot)\big) = 
[x_0; y_0, y_1]g\, [x_0,x_1]y + [x_0, x_1; y_1]g, \]
and when $n=2$,
\begin{align*}
[x_0,x_1,x_2]g\big(\cdot,y(\cdot)\big) = & \  [x_0;y_0,y_2]g\, [x_0,x_1,x_2]y \\
& + [x_0;y_0,y_1,y_2]g\,  [x_0,x_1]y\, [x_1,x_2]y  \\
& + [x_0, x_1; y_1, y_2]g\, [x_1,x_2]y             \\
& + [x_0, x_1, x_2; y_2]g.\\
\end{align*}

In case $y$ is implicitly defined by $g\big(x,y(x)\big) = 0$, the left hand side of Equation \ref{eq:CR1} is zero. In the case $n=1$, therefore, we see that 
\begin{equation}\label{eq:RecursionFormulaZ}
[01]y = - \frac{[01;1]g}{[0;01]g},
\end{equation}
where we now used the shorthand notation from Equations \ref{eq:divdiff-short1} and \ref{eq:divdiff-short2}. For $n\geq 2$, the highest order divided difference of $y$ present in the right hand side of Equation \ref{eq:CR1} appears in the term $[0;0n]g\, [01\cdots n]y$. Moving this term to the left hand side and dividing by $-[0; 0n]g$, one finds a formula that expresses $[01\cdots n]y$ recursively in terms of lower order divided differences of $y$ and divided differences of $g$,
\begin{equation}\label{eq:RecursionFormulaA}
[01\cdots n]y =  - \sum_{k = 2}^n\ \sum_{0 = i_0 < \cdots < i_k = n}\ \sum_{\substack{s = 0\\s = i_s - i_0}}^k
\end{equation}
\[ \hfill \frac{[01\cdots s;i_s i_{s+1}\cdots i_k]g}{[0;0n]g}\!\prod_{l = s + 1}^k [i_{l-1} (i_{l-1} + 1)\cdots i_l]y. \]

We shall now simplify Equation \ref{eq:RecursionFormulaA}. By Equation \ref{eq:RecursionFormulaZ}, the first order divided differences of $y$ appearing in the product of Equation \ref{eq:RecursionFormulaA} can be expressed as quotients of divided differences of $g$. To separate, for every sequence $(i_0,i_1,\ldots,i_k)$ appearing in Equation \ref{eq:RecursionFormulaA}, the divided differences of $g$ from those of $y$, we define an expression involving only divided differences of $g$,
\begin{equation}\label{eq:curlybrackets}
\{i_0\cdots i_k\}g :=  -\sum_{\substack{s = 0 \\s = i_s - i_0}}^k \frac{[i_0\cdots i_s;i_s\cdots i_k]g}{[i_0; i_0 i_k]g} \prod_{\substack{l = s + 1\\ i_l - i_{l-1} = 1}}^k \left( - \frac{ [i_{l-1} i_l; i_l]g}{[i_{l-1}; i_{l-1} i_l]g} \right).
\end{equation}
Note that if a sequence $(i_0,\ldots,i_k)$ starts with precisely $s$ consecutive integers, the expression $\{i_0\cdots i_k\}g$ will comprise $s$ terms.
For instance, 
\begin{align*}
\{023\}g = & \frac{[0;023]g}{[0;03]g}\fg{2}{3}, \\
\{013\}g = & \frac{[0;013]g}{[0;03]g}\fg{0}{1} - \frac{[01;13]g}{[0;03]g}, \\
\{012\}g = &  - \frac{[0;012]g}{[0;02]g}\fg{0}{1}\fg{1}{2} + \frac{[01;12]g}{[0;02]g}\fg{1}{2} - \frac{[012;2]g}{[0;02]g} .
\end{align*}

The remaining divided differences $[i_{l-1}\cdots i_l]y$ in the product of Equation \ref{eq:RecursionFormulaA} are those with $i_l - i_{l-1}\geq 2$, and each of these comes after any $s$ satisfying $s = i_s - i_0$. We might therefore as well start the product of these remaining divided differences at $l=1$ instead of at $l = s + 1$, which has the advantage of making it independent of $s$. Equation \ref{eq:RecursionFormulaA} can thus be rewritten as
\begin{equation}\tag{\ref{eq:RecursionFormulaA}\cprime}\label{eq:RecursionFormulaB}
[0\cdots n]y = \sum_{k = 2}^n \sum_{0 = i_0 < \cdots < i_k = n} \{i_0\cdots i_k\}g \prod_{\substack{l=1\\ i_l - i_{l-1} \geq 2}}^k [i_{l-1}\cdots i_l]y.
\end{equation}
For $n = 2, 3, 4$ this expression amounts to
\begin{align}
  [012]y = &\, \{012\}g,\label{eq:rec2}\\
 [0123]y = &\, \{0123\}g + \{023\}g\, [012]y + \{013\}g\, [123]y,\label{eq:rec3}\\
[01234]y = &\, \{01234\}g + \{0134\}g\, [123]y + \{034\}g\, [0123]y \label{eq:rec4}\\
         + &\, \{0124\}g\, [234]y + \{0234\}g\, [012]y + \{014\}g\, [1234]y\notag\\
         + &\, \{024\}g\, [012]y\, [234]y. \notag
\end{align}

\section{A Formula for Divided Differences of Implicit Functions}\label{sec:formula}

In this section we shall solve the recursive formula from Equation \ref{eq:RecursionFormulaB}. Repeatedly applying Equation \ref{eq:RecursionFormulaB} to itself yields
\begin{align}\label{eq:MainExamples1}
  [012]y =\,& \{012\}g,\\\label{eq:MainExamples2}
 [0123]y =\,& \{0123\}g + \{023\}g\,\{012\}g + \{013\}g\,\{123\}g,\\\label{eq:MainExamples3}
[01234]y =\,& \{01234\}g + \{0134\}g\,\{123\}g + \{034\}g\,\{013\}g\,\{123\}g \\\notag
         +\,& \{034\}g\,\{0123\}g + \{034\}g\,\{023\}g\,\{012\}g + \{0124\}g\,\{234\}g  \\\notag
         +\,& \{0234\}g\,\{012\}g + \{014\}g\,\{134\}g\,\{123\}g + \{014\}g\,\{1234\}g  \\\notag
         +\,& \{014\}g\,\{124\}g\,\{234\}g + \{024\}g\,\{012\}g\,\{234\}g.
\end{align}
Examining these examples, one finds that each term in the right hand sides of the above formulas corresponds to a partition of a convex polygon in a manner we shall now make precise.

With a sequence of labels $0,1,\ldots,n$ we associate the ordered vertices of a convex polygon. A \emph{partition of a convex polygon} is the result of connecting any pairs of nonadjacent vertices with straight lines, none of which intersect. We refer to these straight lines as the \emph{inner edges} of the partition. We denote the set of all such partitions of a polygon with vertices $0,1,\ldots,n$ by $\PPP(0,1,\ldots,n)$. Every partition $\pi\in \PPP(0,1,\ldots,n)$ is described by its set $F(\pi)$ of (oriented) faces. Each face $f\in F(\pi)$ is defined by some increasing sequence of vertices $i_0,i_1,\ldots,i_k$ of the polygon, i.e., $f = (i_0,i_1,\ldots,i_k)$. We denote the set of edges in $\pi$ by $E(\pi)$.

Let $y$ be a function implicitly defined by $g\big(x,y(x)\big) = 0$ and $(x_0,y_0),\ldots,$ $(x_n,y_n)$ be as in Section \ref{sec:Introduction}.
Equations \ref{eq:MainExamples1}--\ref{eq:MainExamples3} suggest the following theorem.

\begin{figure}
\includegraphics[scale=0.5]{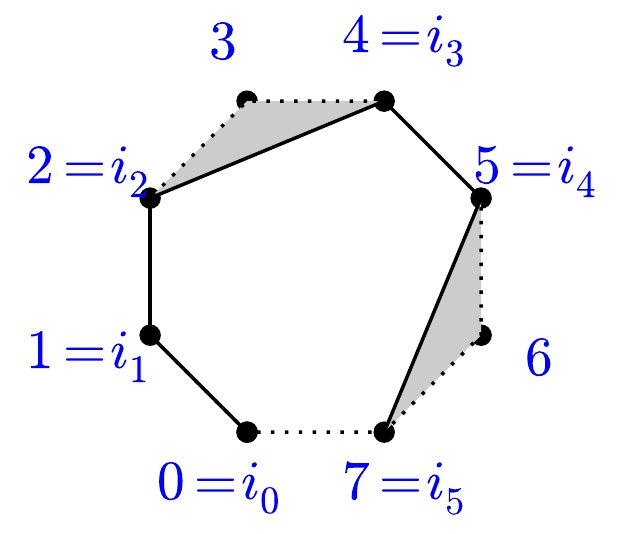}
\caption{For $n=7$, the sequence $\ii = (0,1,2,4,5,7)$ gives rise to the two outer faces $(2,3,4)$ and $(5,6,7)$, which are drawn shaded in the figure. The set $\PPP_\ii$ contains in this case just $1\times 1 = 1$ partition, namely the union of the unique partitions $\{(2,3,4)\}$ and $\{(5,6,7)\}$ of the outer faces.}
\label{fig:inner-outer}
\end{figure}

\begin{theorem}[Main Theorem]\label{thm:main} For $y$ and $g$ defined as above and sufficiently smooth and for $n\geq 2$,
\begin{equation}\label{eq:main}
[0\cdots n]y = \sum_{\pi\in \PPP(0,\ldots,n)}\ \prod_{(v_0,\ldots,v_r)\in F(\pi)} \{v_0\cdots v_r\}g,
\end{equation}
where $\{v_0\cdots v_r\}g$ is defined by Equation \ref{eq:curlybrackets}.
\end{theorem}
Before we proceed with the proof of this theorem, we make some remarks. For $n = 2,3,4$ this theorem reduces to the statements of Equations \mbox{\ref{eq:MainExamples1}--\nolinebreak\ref{eq:MainExamples3}.} To prove Theorem \ref{thm:main}, our plan is to use Equation \ref{eq:RecursionFormulaB} recursively to express $[01\cdots n]y$ solely in terms of divided differences of $g$. We have found it helpful to assign some visual meaning to Equation \ref{eq:RecursionFormulaB}. Every sequence $\ii = (i_0,i_1,\ldots,i_k)$ that appears in Equation \ref{eq:RecursionFormulaB} induces a partition $\pi_\ii \in \PPP(0,1,\ldots,n)$ whose set of faces comprises an \emph{inner face} $(i_0,i_1,\ldots,i_k)$ and \emph{outer faces} $(i_j, i_j + 1, \ldots, i_{j+1})$ for every $j = 0,\ldots, k-1$ with $i_{j+1} - i_j \geq 2$. We denote by $\PPP_{\ii}$ the set of all partitions of the disjoint union of these outer faces. An example of such a sequence $\ii$, together with its inner face, outer faces, and partition set $\PPP_{\ii}$ is given in Figure \ref{fig:inner-outer}.

We shall now associate divided differences to these geometric objects. To each outer face $(i_j, i_j + 1, \ldots, i_{j+1})$ we associate the divided difference $[i_j (i_j + 1) \cdots i_{j+1}]y$, and to each inner face $(i_0,i_1,\ldots,i_k)$ we associate the expression $\{i_0\cdots i_k\}g$. For any sequence $\ii$ that appears in the sum of Equation \ref{eq:RecursionFormulaB}, the corresponding inner face therefore represents that part of Equation \ref{eq:RecursionFormulaB} that can be written solely in terms of divided differences of $g$, while the outer faces represent the part that is still expressed as a divided difference of $y$. 

\begin{figure}
\includegraphics[scale=0.61]{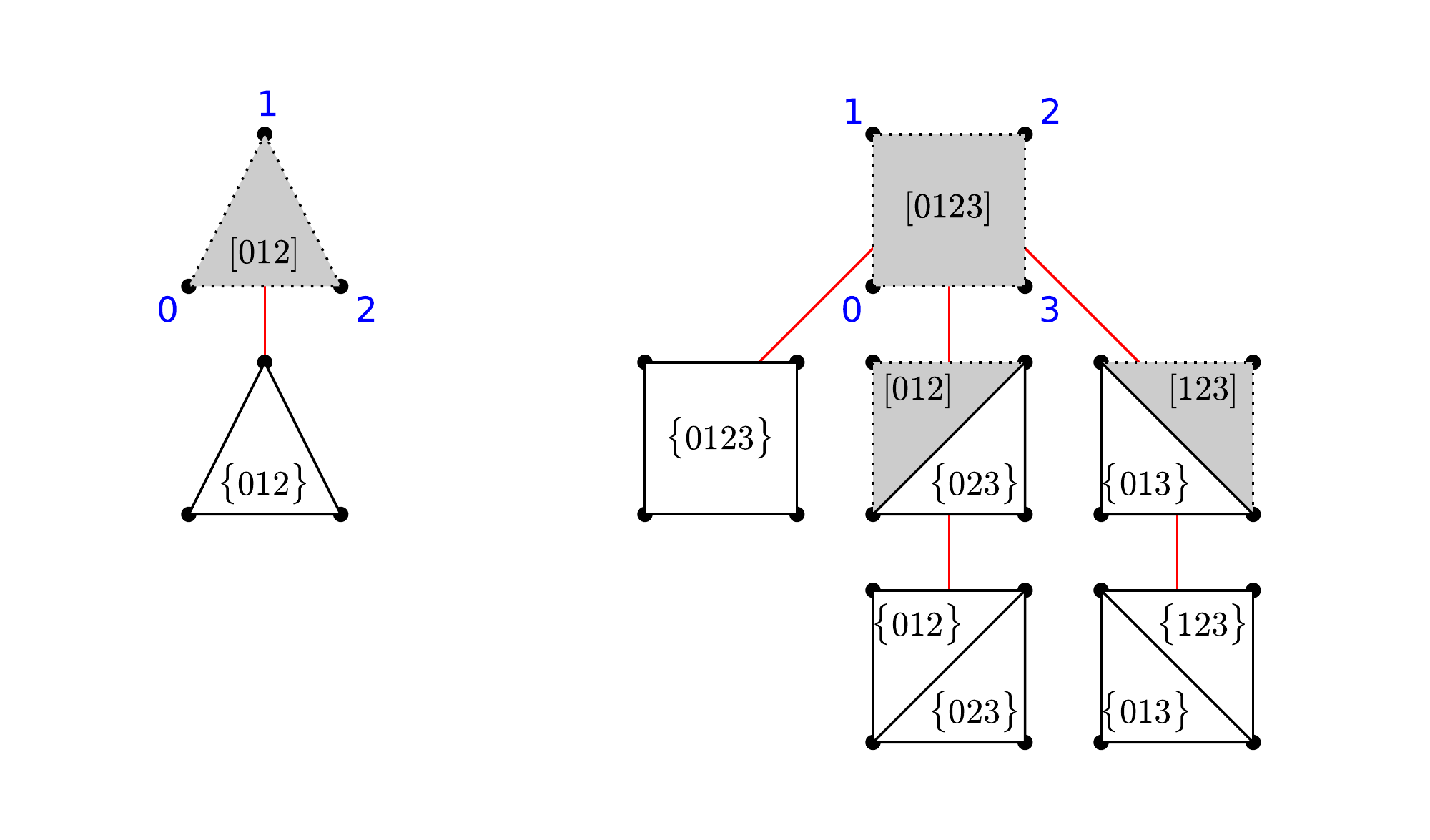}
\hspace*{-3.5em}\includegraphics[scale=0.61]{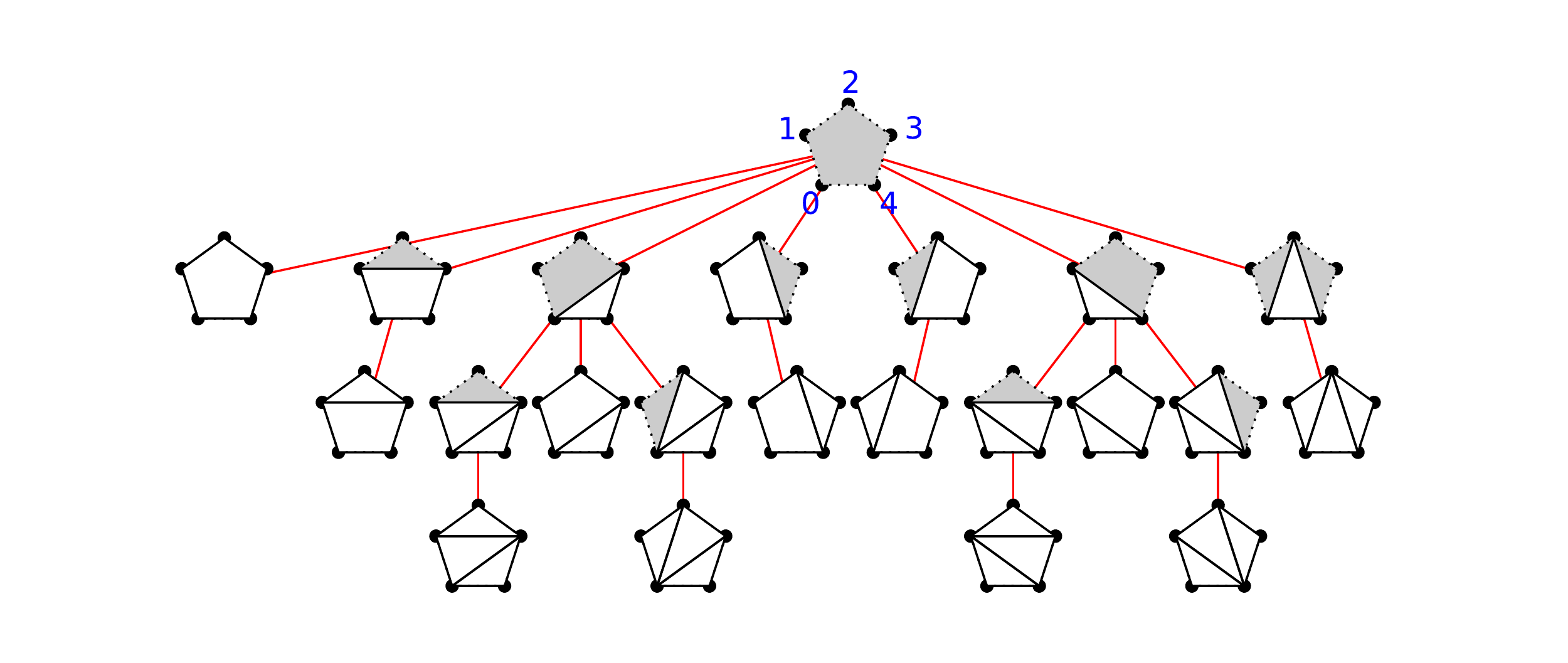}
\caption{For $n=2,3,4$, the figure depicts the recursion trees obtained by repeatedly applying Equation \protect\ref{eq:RecursionFormulaB}. The top levels of these recursion trees correspond to Equations \protect\ref{eq:rec2}--\protect\ref{eq:rec4}.}
\label{fig:RecursionTrees}
\end{figure}

Repeatedly applying Equation \ref{eq:RecursionFormulaB} yields a recursion tree, in which each node represents a product of divided difference expressions associated to inner and outer faces. These recursion trees are depicted in Figure \ref{fig:RecursionTrees} for $n=2,3,4$. Equation \ref{eq:RecursionFormulaB} roughly states that the expression of any nonleaf vertex is equal to the sum of the expressions of its descendants.

\begin{proof}[Proof of the Main Theorem]
This theorem is a generalization of Theorem 1 in \cite{DividedDiffInverse}, and the proofs are analogous. We prove the formula by induction on the order $n$ of the divided difference of $y$.

By the above discussion, the formula holds for $n=2,3,4$. For $n\geq 5$, assume the formula holds for all smaller $n$. Consider the recursive formula from Equation \ref{eq:RecursionFormulaB}. For every sequence $\ii$ that appears in this equation, the corresponding outer faces have fewer vertices than the full polygon. By the induction hypothesis, we can therefore replace each divided difference $[i_l\cdots i_{l+1}]y$ appearing in the product of Equation \ref{eq:RecursionFormulaB} by an expression involving only divided differences of $g$.

As before, let $\PPP_{\ii}$ denote the set of all partitions of the disjoint union of the outer faces induced by $\ii$. Then, by the induction hypothesis, the product in Equation \ref{eq:RecursionFormulaB} is equal to
\[ \sum_{\pi\in \PPP_\ii} \prod_{(v_0,\ldots,v_r)\in F(\pi)} \{v_0\cdots v_r\}g. \]
For a given inner face $\ii$, the set $\PPP_\ii$ can be identified with $\{\pi \in \PPP(0,\ldots,n) : \ii \in F(\pi)\}$ by the bijection $F(\pi)\mapsto F(\pi) \cup \{\ii\}$. Substituting the above expression into Equation \ref{eq:RecursionFormulaB} then yields
\begin{align*} [0\cdots n]y
= & \sum_{\substack{\textup{inner faces}\\\ii=(i_0,\ldots,i_k)}} \{i_0\cdots i_k\}g \sum_{\pi\in \PPP_\ii}\ \prod_{(v_0,\ldots,v_r)\in F(\pi)} \{v_0\cdots v_r\}g \\
= & \sum_{\substack{\textup{inner faces}\\\ii=(i_0,\ldots,i_k)}}\ \sum_{\substack{\pi\in \PPP(0,\ldots,n)\\\ii\in F(\pi)}}\ \prod_{(v_0,\ldots,v_r)\in F(\pi)} \{v_0\cdots v_r\}g \\
= & \sum_{\pi\in \PPP(0,\ldots,n)}\ \prod_{(v_0,\ldots,v_r)\in F(\pi)} \{v_0\cdots v_r\}g.\qedhere
\end{align*}
\end{proof}

Intuitively, this proof can be expressed in terms of the recursion tree as follows. As remarked in the previous section, Equation \ref{eq:RecursionFormulaB} states that the expression of any nonleaf vertex is equal to the sum of the expressions of its descendants. By induction, the expression $[01\cdots n]y$ of the root vertex is therefore equal to the sum of the expressions of the leaves, which, by construction, correspond to partitions of the full polygon. 

\newpage
\begin{example}\label{ex:MainTheorem} Let us apply Theorem \ref{thm:main} to find a simple expression for divided differences of the function $y(x) = \sqrt{1-x^2}$ defined on the interval $(-1,1)$. This function is implicitly defined by the polynomial $g(x,y) = x^2 + y^2 - 1 = 0$. For any knots $x_a, x_b, x_c, x_d$ satisfying $-1 < x_a \leq x_b \leq x_c \leq x_d < 1$ and corresponding function values $y_a,y_b,y_c,y_d$, one finds
\[ [x_a, x_b; y_c]g = x_a + x_b, \qquad   [x_a, x_b, x_c; y_d]g = 1, \]
\[ [x_a; y_b, y_c]g = y_b + y_c, \qquad \ [x_a; y_b, y_c, y_d]g = 1, \]
and all other divided differences of $g$ of nonzero order are zero. In particular, every divided difference of $g$ of total order at least three is zero, which means that the sum in Equation \ref{eq:main} will only be over \emph{triangulations} (i.e., partitions in which all faces are triangles). For a polygon with vertices $0,1,\ldots,n$, Exercise 6.19a of \cite{Stanley2} states that the number of such triangulations is given by the Catalan number
\[ C(n-1) =\frac{1}{n}{2n - 2\choose n - 1}.\]

Consider, for a given triangulation $\pi\in \PPP(0,1,\ldots,n)$, a face $(a,b,c)\in F(\pi)$ from the product in Equation \ref{eq:main}. As any divided difference of the form $[x_a,x_b;y_b,y_c]g$ is zero for this $g$, Equation \ref{eq:curlybrackets} expresses $\{abc\}g$ as a sum of at most two terms. There are four cases.
\[ \{abc\}g = \left\{ \begin{array}{ll}
\displaystyle \frac{-1}{y_a + y_c} \left[ 1 + \frac{x_a + x_b}{y_a + y_b}\cdot \frac{x_b + x_c}{y_b + y_c} \right] & a,b,c~\text{consecutive};\\
\displaystyle \frac{ 1}{y_a + y_c}\cdot\frac{x_a + x_b}{y_a + y_b} & \text{only}~a,b~\text{consecutive};\\
\displaystyle \frac{ 1}{y_a + y_c}\cdot\frac{x_b + x_c}{y_b + y_c} & \text{only}~b,c~\text{consecutive};\\
\displaystyle \frac{-1}{y_a + y_c} & \text{otherwise}.
\end{array}\right. \]
For example, when $n=3$, our convex polygon is a quadrilateral, which admits $C(3-1) = 2$ triangulations $\pi_1$ and $\pi_2$ with sets of faces
\[ F(\pi_1) = \{ (0,1,2), (0,2,3) \},\qquad F(\pi_2) = \{ (0,1,3), (1,2,3) \}. \]
One finds
\[ [x_0,x_1,x_2,x_3]\sqrt{1-x^2} = \{012\}g\, \{023\}g + \{013\}g\,\{123\}g =  \]
\[ \hfill \frac{-1}{(y_0 + y_3)(y_0 + y_2)} \left[ 1 + \frac{x_0 + x_1}{y_0 + y_1}\cdot \frac{x_1 + x_2}{y_1 + y_2} \right] \cdot \frac{x_2 + x_3}{y_2 + y_3} + \]
\[ \hfill \frac{-1}{(y_0 + y_3)(y_1 + y_3)} \left[ 1 + \frac{x_1 + x_2}{y_1 + y_2}\cdot \frac{x_2 + x_3}{y_2 + y_3} \right] \cdot \frac{x_0 + x_1}{y_0 + y_1}.\ \!\  \]
\end{example}

\begin{example}
Next we show that Theorem \ref{thm:main} is a generalization of Theorem 1 of \cite{DividedDiffInverse}, which gives a similar formula for inverse functions. To see this, we apply Theorem \ref{thm:main} to a function $y$ implicitly defined by a function $g(x,y) = x - h(y)$. Referring to Equation \ref{eq:curlybrackets}, we need to compute $[i_0\cdots i_s;i_s\cdots i_k]g$ for this choice of $g$ and various indices $i_0,\ldots,i_k$ and $s\in \{0,\ldots,k\}$. Applying the recursive definition of bivariate divided differences, one obtains
\begin{align*}
[i_0\cdots i_s;i_s\cdots i_k]x \hfill = & 
\left\{ \begin{array}{ll} x_{i_0} & \textup{if}~s = 0, s = k;\\ 1 & \textup{if}~s = 1, s = k;\\  0 & \textup{otherwise}, \end{array} \right. \\
 [i_0\cdots i_s;i_s\cdots i_k]h(y) = &
\left\{ \begin{array}{ll} [i_s\cdots i_k]h & \textup{if}~s = 0;\\ 0 & \textup{otherwise}. \end{array} \right.
\end{align*}

Consider a face $f = (v_0,\ldots,v_r)$ of a given partition $\pi \in \PPP(0,\ldots,n)$ in Equation \ref{eq:main}. Since $r\geq 2$, the divided difference $[v_0\cdots v_s;v_s\cdots v_r]\big(x-h(y)\big)$ is zero for $s\geq 1$. Using this, Equation \ref{eq:curlybrackets} expresses $\{v_0\cdots v_r\}g$ as a single term 
\begin{align*}
\{v_0\cdots v_r\}g
 = & -\frac{[v_0;v_0\cdots v_r]g}{[v_0;v_0v_r]g} \prod_{\substack{l=1\\v_l - v_{l-1} = 1}}^r \left(-\frac{[v_{l-1}v_l;v_l]g}{[v_{l-1};v_{l-1}v_l]g} \right) \\
 = & -\frac{[v_0\cdots v_r]h}{[v_0v_r]h}         \prod_{\substack{l=1\\v_l - v_{l-1} = 1}}^r \frac{1}{[v_{l-1}v_l]h}.
\end{align*}
Taking the product over all faces in the partition $\pi$, the denominators of the factors in the above equation correspond to the edges of the partition, while the numerators correspond to the faces of the partition. As there is a minus sign for each face in the partition, we arrive at the formula 
\begin{equation}\label{eq:inverse}
[01\cdots n]y = \sum_{\pi\in \PPP(0,\ldots,n)} (-1)^{\#F(\pi)}
\frac{\displaystyle \prod_{(v_0,v_1,\ldots,v_r)\in F(\pi)} [v_0v_1\cdots v_r]h}{\displaystyle \prod_{(v_0,v_1)\in E(\pi)} [v_0 v_1]h},
\end{equation}
which appears as Equation 11 in \cite{DividedDiffInverse}.
\end{example}

Note that the inverse of the algebraic function $y = \sqrt{1-x^2}$ in Example \ref{ex:MainTheorem} is again an algebraic function. Equation \ref{eq:inverse} would therefore not have been of much help to find a simple expression for divided differences of $y$. In fact, Example \ref{ex:MainTheorem} can be thought of as one of the simplest examples for which Theorem \ref{thm:main} improves on Equation \ref{eq:inverse}, as it concerns a polynomial $g$ with bidegree as low as (2,2).

\begin{example} In this example we shall derive a quotient rule for divided differences. That is, we shall find a formula that expresses divided differences of the quotient $y = P(x)/Q(x)$ in terms of divided differences of $P$ and of $Q$. Let $g(x,y) = Q(x)y - P(x)$. Then, in Equation \ref{eq:curlybrackets},
\begin{equation}\label{eq:quotient1}
[i_0\cdots i_s;i_s\cdots i_k]g = 
\left\{ \begin{array}{rl}
y_{i_s}[i_0\cdots i_s]Q - [i_0\cdots i_s]P & \textup{if}~s = k;\\
\,[i_0\cdots i_s]Q & \textup{if}~s = k-1;\\
0 & \textup{otherwise}.
\end{array} \right.
\end{equation}

In Equation \ref{eq:main}, therefore, the only partitions with a nonzero contribution are those whose faces have all their vertices consecutive, except possibly the final one. In particular, the inner face with vertices $0 = i_0 < \cdots < i_k = n$ should either be the full polygon, or should have a unique inner edge $(i_{k-1}, n)$. By induction, it follows that the partitions with a nonzero contribution to Equation \ref{eq:main} are precisely those for which all inner edges end at $n$. These partitions correspond to subsets $I \subset \{1,2,\ldots,n-2\}$, including the empty set, by associating with any such $I$ the partition with inner edges $\{ (i,n)\,:i\, \in I\}$.
Equation \ref{eq:main} becomes
\begin{align}\label{eq:quotient2}
[0\cdots n]\frac{P}{Q} =\, & \{0\cdots n\}g + \\\notag
              \, & \sum_{r=1}^{n-2}\ \sum_{k=1}^r\ \sum_{0 = i_0 < i_1 < \cdots < i_k = r} \{r\cdots n\}g\ \prod_{j=1}^k \{i_{j-1}\cdots i_j n\}g,
\end{align}
where the dots represent consecutive nodes and an empty product is understood to be one. A long but straightforward calculation involving Equations \ref{eq:curlybrackets}, \ref{eq:quotient1}, and \ref{eq:quotient2} yields
\begin{align*}
[0\cdots n] \frac{P}{Q} =\, & \frac{[0\cdots n]P}{Q_0} + \\
 \, & \sum_{r=1}^n \frac{[r\cdots n]P}{Q_r}  \sum_{k=1}^r (-1)^k \!\!\sum_{0 = i_0 < i_1 < \cdots < i_k = r}\ \prod_{j=1}^k \frac{[i_{j-1}\cdots i_j]Q}{Q_{i_{j-1}}},
\end{align*}
where $Q_i := Q(x_i)$ for $i = 0,\ldots, n$. Alternatively, this equation can be found by applying a univariate chain rule to the composition $x\mapsto Q(x)\mapsto 1/Q(x)$, as described in Section 4 of \cite{DividedDiffChain}.
\end{example}

Finally, we note that taking the limit $x_0,\ldots,x_n\to x$ in  Equations \ref{eq:RecursionFormulaZ}, \ref{eq:curlybrackets}, \ref{eq:MainExamples1}, and \ref{eq:MainExamples2} yields
\begin{align*}
  y'(x)   = & - \frac{g_{10}}{g_{01}},\\
  y''(x)  = & - \frac{g_{20}}{g_{01}} + 2\frac{g_{11} g_{10}}{g_{01}^2} - \frac{g_{02} g_{10}^2}{g_{01}^3},\\
  y'''(x) = & - \frac{g_{30}}{g_{01}} + 3\frac{g_{21} g_{10}}{g_{01}^2} + 3\frac{g_{20} g_{11}}{g_{01}^2}
              -3\frac{g_{20} g_{10} g_{02}}{g_{01}^3} -3\frac{g_{12} g_{10}^2 }{g_{01}^3} \\
            & -6\frac{g_{11}^2 g_{10}}{g_{01}^3} + \frac{g_{10}^3 g_{03}}{g_{01}^4} + 9\frac{g_{11} g_{10}^2 g_{02}}{g_{01}^4}
              -3\frac{g_{10}^3 g_{02}^2}{g_{01}^5},\\
\end{align*}
where we introduced the shorthand
\[ g_{st} := \frac{\partial^{s+t} g}{\partial x^s \partial y^t}\big(x,y(x)\big).\]
These formulas agree with the examples given in \cite{ComtetFiolet}, \cite[Page 153]{Comtet} and with a formula stated as Equation 7 in \cite{Wilde}.

\section*{Acknowledgment}
We wish to thank Paul Kettler, whose keen eye for detail provided us with many valuable comments on a draft of this paper.

\section*{References}

\bibliographystyle{amsxport}

\end{document}